\documentclass[12pt,leqno,draft]{article}
\usepackage{amsfonts}
\usepackage{amsmath, amsthm, amsfonts, amssymb, color}
\usepackage{mathrsfs}
\usepackage{color}
\theoremstyle{definition}

\newcommand{\scr}[1]{\mathscr #1}
\definecolor{wco}{rgb}{0.5,0.2,0.3}

\numberwithin{equation}{section} \theoremstyle{remark}

\newcommand{\ua}{\uparrow}

\title{{\bf 
 Super and Weak Poincar\'e Inequalities for Sticky-Reflected Diffusion Processes}\footnote{Supported in
 part by  National Key R\&D Program of China (No. 2022YFA1006000) and State Key Lab. Synthetic Biology.} }
\author{
{\bf    Feng-Yu Wang  }\\
\footnotesize{ Center for Applied Mathematics and KL-AAGDM, Tianjin
University, Tianjin 300072, China} \\
\footnotesize{  wangfy@tju.edu.cn}\\
}
\begin{document}
\allowdisplaybreaks
\def\R{\mathbb R}  \def\ff{\frac} \def\ss{\sqrt} \def\B{\mathbf
B} \def\W{\mathbb W}
\def\N{\mathbb N} \def\kk{\kappa} \def\m{{\bf m}}
 \def\ddd{D^*}
\def\dd{\delta} \def\DD{\Delta} \def\vv{\varepsilon} \def\rr{\rho}
\def\<{\langle} \def\>{\rangle} \def\GG{\Gamma} \def\gg{\gamma}
  \def\nn{\nabla} \def\pp{\partial} \def\E{\mathbb E}
\def\d{\text{\rm{d}}} \def\bb{\beta}  \def\D{\scr D}
  \def\si{\sigma} \def\ess{\text{\rm{ess}}}
\def\beg{\begin} \def\beq{\begin{equation}}  \def\F{\scr F}
\def\Ric{\text{\rm{Ric}}} \def\Hess{\text{\rm{Hess}}}
\def\e{\text{\rm{e}}} \def\ua{\underline a} \def\OO{\Omega}  \def\oo{\omega}
 \def\tt{\tilde} \def\Ric{\text{\rm{Ric}}}
\def\cut{\text{\rm{cut}}} \def\P{\mathbb P} \def\ifn{I_n(f^{\bigotimes n})}
\def\C{\scr C}      \def\thetaa{\mathbf{r}}     \def\r{r}
\def\gap{\text{\rm{gap}}} \def\prr{\pi_{{\bf m},\varrho}}  \def\r{\mathbf r}
\def\Z{\mathbb Z} \def\vrr{\varrho} \def\ll{\lambda}
\def\L{\scr L}\def\Tt{\tt} \def\TT{\tt}\def\II{\mathbb I}
\def\i{{\rm in}}\def\Sect{{\rm Sect}}  \def\H{\mathbb H}
\def\M{\scr M}\def\Q{\mathbb Q} \def\texto{\text{o}} \def\LL{\Lambda}
\def\Rank{{\rm Rank}} \def\B{\scr B} \def\i{{\rm i}} \def\HR{\hat{\R}^d}
\def\to{\rightarrow}\def\l{\ell}\def\iint{\int}
\def\EE{\scr E}\def\Cut{{\rm Cut}}
\def\A{\scr A} \def\Lip{{\rm Lip}}
\def\BB{\scr B}\def\Ent{{\rm Ent}}\def\L{\scr L}
\def\R{\mathbb R}  \def\ff{\frac} \def\ss{\sqrt} \def\B{\mathbf
B}
\def\N{\mathbb N} \def\kk{\kappa} \def\m{{\bf m}}
\def\dd{\delta} \def\DD{\Delta} \def\vv{\varepsilon} \def\rr{\rho}
\def\<{\langle} \def\>{\rangle} \def\GG{\Gamma} \def\gg{\gamma}
  \def\nn{\nabla} \def\pp{\partial} \def\E{\mathbb E}
\def\d{\text{\rm{d}}} \def\bb{\beta}  \def\D{\scr D}
  \def\si{\sigma} \def\ess{\text{\rm{ess}}}
\def\beg{\begin} \def\beq{\begin{equation}}  \def\F{\scr F}
\def\Ric{\text{\rm{Ric}}} \def\Hess{\text{\rm{Hess}}}
\def\e{\text{\rm{e}}} \def\ua{\underline a} \def\OO{\Omega}  \def\oo{\omega}
 \def\tt{\tilde} \def\Ric{\text{\rm{Ric}}}
\def\cut{\text{\rm{cut}}} \def\P{\mathbb P} \def\ifn{I_n(f^{\bigotimes n})}
\def\C{\scr C}      \def\thetaa{\mathbf{r}}     \def\r{r}
\def\gap{\text{\rm{gap}}} \def\prr{\pi_{{\bf m},\varrho}}  \def\r{\mathbf r}
\def\Z{\mathbb Z} \def\vrr{\varrho} \def\ll{\lambda}
\def\L{\scr L}\def\Tt{\tt} \def\TT{\tt}\def\II{\mathbb I}
\def\i{{\rm in}}\def\Sect{{\rm Sect}}  \def\H{\mathbb H}
\def\M{\scr M}\def\Q{\mathbb Q} \def\texto{\text{o}} \def\LL{\Lambda}
\def\Rank{{\rm Rank}} \def\B{\scr B} \def\i{{\rm i}} \def\HR{\hat{\R}^d}
\def\to{\rightarrow}\def\l{\ell}\def\BB{\mathbb B}
\def\8{\infty}\def\I{1}\def\U{\scr U} \def\n{{\mathbf n}}\def\v{V}\def\aa{\alpha}
\maketitle

\begin{abstract} As a continuation to \cite{MRW} where the Poincar\'e and log-Sobolev inequalities were studied for the sticky-reflected Brownian motion on Riemannian manifolds with boundary,
this paper establishes the  super and weak Poincar\'e inequalities for more general sticky-reflected diffusion processes. As applications, the convergence  rate and uniform integrability of the associated diffusion semigroups are characterized.
The main results are illustrated by concrete examples.
 \end{abstract} \noindent
 \vskip 2cm

\section{Introduction}

 Let $(M,\<\cdot,\cdot\>)$ be a $d$-dimensional open Riemannian manifold with smooth boundary $(\pp M,\<\cdot,\cdot\>_\pp)$. We consider a Markov process on $\bar M$ as follows:
 \beg{enumerate}\item[(1)] Starting from a point in $M$ it moves as a diffusion process in $M$ until hits the boundary $\pp M$; if the staring point is on the boundary, the hitting time is $0$.
 \item[$(2)$] From the hitting time to $\pp M$, it stays at the hitting point (i.e. without boundary diffusion), or it moves as another diffusion process on $\pp M$ (i.e. with boundary diffusion), until a random time determined by
 the strength of reflection.
   \item[$(3)$] At the random time,  it is reflected into $M$ and moves as the diffusion in $M$ again until hits the boundary, and repeatedly.
   \end{enumerate}

  This process is called a sticky-reflected diffusion process, or diffusion process with Wentzell's boundary condition since the study goes back to   Wentzell  \cite{Wen},  and   has been used to describe   interacting particle systems with singular boundary or zero-range pair interaction, see \cite{AD, DG, GD, KR} and references therein. Rigorous constructions of  sticky-reflected diffusion processes  were presented in \cite{I, IW, TW, W} for $M$ being a special domain (e.g. ball) in $\R^d$, and in \cite{CP,T} for more general domains. See  \cite{FG,GV} for the study by using Dirichlet forms.

Recently, optimal constants in the  Poincar\'e and log-Sobolev inequalities have been estimated in \cite{M,MRW} for the sticky-reflected  (weighted) Brownian motions on $\bar M$, which extend the  corresponding results derived in  \cite{KMV} for strictly convex manifolds with positive curvature.  In this paper, we study the super and weak Poincar\'e inequalities, which were introduced in \cite{W00} and \cite{RW01} respectively,
for the sticky-reflected diffusion processes.

 Let $\LL$ and $\LL_\pp$ be   the volume measures on $\bar M$ and $\pp M$ respectively.  Let $V\in C^1(M)$ and $W\in C^1(\pp M)$ be   such that
 $$Z_V:= \int_M \e^{V}\d\LL<\infty,\ \ \ Z_W^\pp :=\int_{\pp M} \e^W \d\LL_\pp<\infty.$$
 Then
$$\mu_V(\d x):=\ff 1{Z_V}\e^{V(x)}\LL(\d x),\ \ \ \mu_W^\pp (\d x):= \ff {1_{\pp M}(x)} {Z_W^\pp}\e^{W(x)}\LL_\pp(\d x)$$
 are probability measures on $\bar M$, where $\mu_W^\pp$ is fully supported on the boundary  $\pp M$.

 For two constants $\dd\in [0,\infty),\gg\in (0,\infty)$,  we consider the operator
  $$\scr L:= 1_{M}(\DD+\nn V)  + 1_{\pp M}\big(\dd [\DD^\pp+\nn^\pp W] + \gg \e^{V-W} N\big),$$
 where $\DD$ and $\nn$ are the Laplacian and gradient operators in $M$,   $\DD^\pp$ and $\nn^\pp$ are the corresponding ones on $\pp M$, and $N$ is the unit inward normal vector field on $\pp M$.
 The diffusion process generated by $\scr L$ is called  sticky-reflected diffusions with inside diffusion generated by $\DD+\nn V$ and boundary diffusion generated by $\dd(\DD^\pp+\nn^\pp W).$
 The constant $\gg>0$ measures the strength of reflection, and the model converges to the reflected diffusion process as   $\gg\to\infty$.
 When $\dd =0,$ the process is called sticky-reflected diffusion process without boundary diffusion, and if moreover $\gg=0$ it becomes the diffusion with absorbing (i.e. killed) boundary.

 To formulate the associated Dirichlet form,
 let
 \beq\label{TH} \theta:=\ff{\gg Z_W^\pp}{\gg Z_W^\pp+ Z_V}\in (0,1).\end{equation}
 Then the associate invariant probability measure for the sticky-reflected diffusion process is the following convex combination of $\mu_V$ and $\mu_W^\pp$:
 $$\mu:= \theta \mu_V+ (1-\theta)\mu_W^\pp.$$
 Indeed,  by the integration by parts formula, for any $f,g\in C_0^2(\bar M)$, the class of
 $C^2$-functions on $\bar M$ with compact support, we have
 $$-\int_{\bar M} (f \scr L g)\d\mu=\EE_\dd(f,g):=\int_{\bar M} \big\{\<\nn f,\nn g\>+ \dd\<\nn^\pp f,\nn^\pp g\>_\pp\big\}\d\mu.$$
 Then it is standard   that the form $(\EE_\dd,C_0^2(\bar M))$ is closable and its closure is a regular symmetric local Dirichlet form in $L^2(\mu)$, so that it associates with a unique diffusion process
 on $\bar M$ with generator $\scr L$, see \cite{FU}.
In particular, when $M$ is a smooth domain in $\R^d$, the sticky-reflected diffusion process can be constructed by solving the SDE
 \beg{align*} &\d X_t=   1_{M}(X_t)\Big\{\ss 2\d B_t +\nn V(X_t)\d t\Big\}\\
&+ 1_{\pp M}(X_t)\Big\{\ss{2\dd} P(X_t)\circ\d B_t+\dd \nn^\pp W(X_t)\d t+ \gg  \e^{V-W}N(X_t)\d t\Big\},\ \ t\ge 0, \end{align*}
 where $B_t$ is the $d$-dimensional Brownian motion, $\d B_t$ and $\circ\d B_t$ are It\^o's and Stratonovich's differentials respectively,   $N$ is the inward unit normal vector field on $\pp M$, and for each $x\in \pp M$,
 $$P(x): \R^d\to T_x\pp M$$ is the orthogonal  projection operator.
 According to \cite{GV}, $\theta$ in \eqref{TH} is the average time for $X_t$ staying in $M$:
 $$\lim_{t\to\infty}\ff 1 t\int_0^t 1_M(X_s)\d s= \theta.$$

  On the other hand, functional inequalities are power  tools  in the study of Markov processes,
for instances, the Sobolev/Nash type inequality characterizes heat kernel estimates (see e.g. \cite{Davies}), Gross' log-Sobolev inequality \cite{Gross1, Gross2} describes the hypercontractivity and exponential ergodicity in entropy, the Poincar\'e (spectral gap) inequality is equivalent to the exponential ergodicity in $L^2$,
the super Poincar\'e inequality introduced in \cite{W00} is equivalent to the empty of essential spectrum and uniform integrability of semigroup, and   the weak Poincar\'e inequality introduced in \cite{RW01} describes general ergodicity rate which is slower than exponential. Recently, by using log-Sobolev inequalities,  L. Gorss \cite{Gross3} showed the invariance of intrinsic ultracontractivity for Schr\"odinger operators under a class of unbounded perturbations, so that   the corresponding result  in \cite{DS} under bounded perturbations is improved. 

In the following, we   simply denote $\nu(f)=\int_{\bar M}f\d\mu$ for a measure
 $\nu$ on $\bar M$ and a function $f\in L^1(\nu).$ Let
$$ |\nn f|^2:=\<\nn f,\nn f\>,\ \  |\nn^\pp f|_\pp^2:=\<\nn^\pp f,\nn^\pp f\>_\pp.$$
We have
$$\EE_\dd(f,f)= \theta\mu_V(|\nn f|^2)+ (1-\theta)\dd \mu_W^\pp(|\nn^\pp f|_{\pp}^2).$$ In this paper, we investigate  the super Poincar\'e inequality
\beq\label{SP} \mu(f^2)\le r\EE_\dd(f,f) +\bb(r)\mu(|f|)^2,\ \ \ r>0,\ f\in C_b^1(\bar M)\end{equation}
and the weak Poincar\'e inequality
\beq\label{WP} \mu(f^2)\le \aa(r)\EE_\dd(f,f) + r \|f\|_\infty^2,\ \ \ r>0,\ f\in C_b^1(\bar M), \ \mu(f)=0,\end{equation}
 where $\bb:  (0,\infty)\to (0,\infty)$   is crucial to estimate  the associated diffusion semigroup and higher order eigenvalues of the generator, see  \cite{W00};
 and $\aa:  (0,\infty)\to (0,\infty)$  corresponds to the convergence rate of the associated Markov semigroup, see \cite{RW01}.

We will estimate the smallest rate functions $\aa$ and $\bb$ for the above introduced sticky-reflected diffusion process:
\beg{align*} &\bb(r):= \sup\big\{\mu(f^2)-r\EE_\dd(f,f): \ f\in C_b^1(\bar M),\ \mu(|f|)=1\big\},\\
&\aa(r):= \sup\big\{(\mu(f^2)-r\|f\|_\infty^2)^+: \ f\in C_b^1(\bar M),\ \mu(f)=0,\ \EE_\dd(f,f)=1\big\},\ \ r>0.\end{align*}
It is easy  to see that $\aa(r)=0$ for $r\ge 1$, and $\bb(r)\ge 1 $ with $\bb(r)=1$ for $r\ge C(P)$, where $C(p)$ is the smallest positive constant such that the Poincar\'e inequality
holds:
$$ \mu(f^2)\le C(p)\EE_\dd(f,f) +\mu(f)^2,\ \ f\in C_b^1(\bar M).$$
So, it suffices to estimate $\aa(r)$ and $\bb(r)$ for small $r>0$, in particular, to characterize the rates of
$\aa(r)\uparrow\infty$ and $\bb(r)\uparrow\infty$ as $r\downarrow 0.$

Since the super/weak Poincar\'e inequalities have been well studied for elliptic diffusions on manifolds with Neumann boundary or without boundary, see \cite{W05}, we will estimate $\aa(r)$ and $\bb(r)$ using the rate functions $\bb_V, \bb^\pp_W, \aa_V$ and $\aa^\pp_W$ in the following functional inequalities:
  \beq\label{SP0} \mu_V(f^2)\le r \mu_V(|\nn f|^2)+\bb_V(r)\mu_V(|f|)^2,\ \ r>0,\ f\in C_b^1(\bar M), \end{equation}
\beq\label{SP1'}\mu_W^\pp(f^2)\le r \mu_W^\pp(|\nn^\pp f|_{\pp}^2)+\bb^\pp_W(r)\mu_W^\pp(|f|)^2,\ \ r>0,\ f\in C_b^1(\pp M), \end{equation}
\beq\label{WP0} \mu_V(f^2)\le \aa_V(r) \mu_V(|\nn f|^2)+ r\|f\|_\infty^2,\ \ r>0,\ f\in C_b^1(\bar M),\ \mu_V(f)=0, \end{equation}
\beq\label{WP1'} \mu_W^\pp(f^2)\le \aa_W(r) \mu_W^\pp(|\nn^\pp f|_{\pp}^2)+r\|f\|_\infty^2,\ \ r>0,\ f\in C_b^1(\pp M),\ \mu_W^\pp(f)=0.\end{equation}

In Section 2, we recall some known results on super and weak Poincar\'e inequalities, which will be applied to the sticky-reflected diffusions.
In Section  3 and Section 4   we establish these inequalities for $ \EE_\dd$  with $\dd>0$ (the case with boundary diffusion), and $\dd=0$ (the case without boundary diffusion), respectively.
Some examples are presented to illustrate our main results.

 \section{A review on super and weak Poincar\'e inequalities}

  Let $(E,\F,\mu)$ be a separable probability space,  let  $(\EE,\D(\EE))$
  be a Dirichlet form on $L^2(\mu)$, let $(L,\D(L))$ and $P_t:=\e^{tL}$ be the associated generator and (sub-) Markov  semigroup.  For any $p,q\in [1,\infty]$, let
  $\|\cdot\|_{p\to q}$ denote the operator norm from $L^p(\mu)$ to $L^q(\mu).$
 In this section, we summarize some results on super and weak Poincar\'e inequalities, where detailed proofs can be found in the book \cite{W05} or \cite{W00, RW01}.

\subsection{Super Poincar\'e  inequality }

We say that
$(\EE,\D(\EE))$ satisfies the  super Poincar\'e inequality, if there   exists a (decreasing) function $\bb: (0,\infty)\to (0,\infty)$   such that
\beq\label{SP1} \mu(f^2)\le r\EE(f,f) +\bb(r)\mu(|f|)^2,\ \ \
r>0, f\in \D(\EE).\end{equation}  This inequality was introduced in \cite{W00} to study the essential spectrum of the generator $L$, and has been further used to estimate the semigroup $P_t$.

 We first introduce the link between   \eqref{SP1} and the uniform integrability of $P_t$.

\beg{thm}[\cite{W05},  Lemma 3.3.5, Theorem 3.3.6]\label{TT1}   The following assertions hold.
\beg{enumerate} \item[$(1)$] The inequality
$(\ref{SP1})$ holds if and only if
$$\mu(|P_tf|^2) \le \e^{-2rt} \mu(f^2) +
\bb(r^{-1}) (1-\e^{-2rt})\mu(|f|)^2,\ \ \ r>0, t\ge 0, f\in
L^2(\mu).$$
\item[$(2)$] Let
$$\GG_t(s)= \inf\{r\ge 0: \bb(1/r)(\e^{2rt}-1)\ge s^2\},\ \ s\ge 0.$$
If $(\ref{SP1})$ holds, then
$$\sup_{\mu(f^2)=1}\mu((P_tf)^21_{\{|P_tf|>r\}})\le
\exp[-2t\GG_t(\vv r)]\big/(1-\vv)^2,\ \ \ r>0, \ \vv\in
(0,1).$$
\item[$(3)$] If  $(\EE,\D(\EE))$ is symmetric and there exists $t>0$ such that
$$\phi_t(s):=\sup\{
\mu((P_tf)^21_{\{|P_tf|>s\}}): \mu(f^2)=1\}\to 0$$ as $s\to\infty,$
then $(\ref{SP1})$ holds with
$$\bb(r)= \ff {r [\phi_t^{-1}(\e^{-2t/r}/2)]^2\e^{2t/r}}{4t},\ \ r>0,$$
where $\phi_t^{-1}(r)=\inf\{s>0: \phi_t(s)\le r\}.$\end{enumerate} \end{thm}

\beg{cor}[\cite{W05}, Corollary 3.3.10]\label{CC1}  Let $\dd\in(0,1]$.    If $(\ref{SP1})$ holds with
$$\bb(r)=\exp[c(1+r^{-1/\dd})],\ \ r>0$$ for some $c>0$, then there exists decreasing $C: (0,\infty)\to (0,\infty)$ such that
\beq\label{3.3.19}\sup_{\mu(f^2)=1}\int_E(P_tf)^2\exp\Big\{C_t\big[\log(1+(P_tf)^2)\big]^\dd
\Big\} \d\mu<\infty,\ \ t\in (0,\infty).\end{equation}  On the other hand, if
$(\EE,\D(\EE))$ is symmetric and $(\ref{3.3.19})$ holds for some
$t>0$ and $C_t>0$, then $(\ref{SP1})$ holds with
$\bb(r)=\exp[c(1+r^{-1/\dd})]$ for some $c>0$. \end{cor}

Next, we consider   three boundedness properties of
$P_t$ by using (\ref{SP1}).

\beg{defn} Let $(E,\F,\mu)$ be a measure space and $P_t$ a
semigroup on $L^2(\mu)$ which is bounded on $L^p(\mu)$ for all
$p\in [1,\infty]$. $P_t$ is called  \emph{hyperbounded}
\index{semigroup!hyperbounded} if $\|P_t\|_{2\to 4}<\infty$ for
some $t>0$; \emph{superbounded} if $\|P_t\|_{2\to 4}<\infty$ for
all $t>0$; and \emph{ultrabounded} if
$\|P_t\|_{1\to\infty}<\infty$ for all $t>0.$
\index{semigroup!ultrabounded}\end{defn}

By Riesz-Thorin interpolation theorem, $P_t$ is hyperbounded if
and only if $\|P_t\|_{p\to q}<\infty $ holds for some constants $t>0$ and    $1<p<q<\infty.$  
It is clear that the
superboundedness implies the hyperboundedness, and  they are implied by the ultraboundedness.

 \beg{thm}[\cite{W05}, Theorems 3.3.13, 3.3.14, 3.3.15] \label{TT2}   We have the following assertions on the hyper/super/ultraboundedness of $P_t$.
 \beg{enumerate}\item[$(1)$]
  If $(\ref{SP1})$ holds with $\bb(r)= \exp[c(1+r^{-1})]$ for
some $c>0,$ then $P_t$ is hyperbounded, and the converse result holds if $(\EE,\D(\EE))$
is symmetric.
\item[$(2)$]
  If $(\ref{SP1})$ holds for some $\bb$  with $\lim_{r\to
0}r\log \bb(r)=0$, then $P_t$ is
superbounded. Conversely, if $P_t$ is superbounded and $(\EE,\D(\EE))$ is
symmetric, then $(\ref{SP1})$ holds for
$$\bb(r):= \inf_{s\le r}\bigg(\ff s {3\e}\land 2\bigg)
\inf_{t>0}\ff {(1+\|P_t\|_{2\to 4}^2)^2}t \exp[6t/s],\ \ r>0,$$
which satisfies  $\lim_{r\to 0}r\log \bb(r)=0.$
\item[$(3)$]   If
$(\ref{SP1})$ holds with $\bb$ satisfying
$$\Psi(t):=\int_t^\infty \ff {\bb^{-1}(r)}r \d r<\infty,\ \ \ t>\inf \bb,$$
then $P_t$ is ultrabounded with
$$ \|P_t\|_{1\to \infty}\le \inf_{\vv\in (0,1)}\max\big\{\vv^{-1}\inf\bb,\
 \Psi^{-1}((1-\vv)t)\big\},\ t>0,$$ where
$\Psi^{-1}(t):=\inf\{r\ge \inf\bb: \Psi(r) \le t\}$.  On the other hand, if $P_t$ is ultra-bounded, then $(\ref{SP1})$
holds for
$$\bb(r)= \inf_{s\le r,t>0} \ff {s\|P_t\|_{1\to
\infty}}{t}\exp\big[t/s -1\big],\ \ r>0.$$
\end{enumerate}\end{thm}

The following is a direct consequence of Theorem \ref{TT2}(3).

\beg{cor}[\cite{W05}, Theorem 3.3.15] \label{CC2}  We have the following correspondence between $\bb$ and $\|P_t\|_{1\to\infty}.$
\beg{enumerate}
\item[$(1)$] Let $\dd>1.\ (\ref{SP1})$ with
$\bb(r)=\exp[c(1+r^{-1/\dd})]$ for some $c>0$ is equivalent to
$$\|P_t\|_{1\to\infty}\le \exp[\ll(1+t^{-1/(\dd-1)})],\ \ \ \ t>0,
$$
for some $\ll>0$.
\item[$(2)$] Let $p>0.\ (\ref{SP1})$ with $\bb(r)=c(1+r^{-p/2})$ for some
$c>0$ is equivalent to
\beq\label{3.3.24} \|P_t\|_{1\to\infty}\le
\ll(1+t^{-p/2})\end{equation} for some $\ll>0$ and all $t>0.$  \end{enumerate} \end{cor}

 \subsection{Weak Poincar\'e inequality}

 In this part, we assume that the Dirichlet form
  $(\EE,\D(\EE))$ is   irreducible and  conservative, i.e. $1\in \D(\EE)$ with $\EE(1,1)=0$, and  $f\in\D(\EE)$ with $\EE(f,f)=0$ implies that $f$ is constant.
In this case, we have
$$\lim_{t\to\infty} \mu(|P_tf-\mu(f)|^2) =0,\ \ f\in L^2(\mu).$$
In the following we introduce the link between the convergence rate  for $\|P_t-\mu\|_{\infty\to 2}\to 0$ as $t\to 0$, and the function $\aa: (0,\infty)\to (0,\infty)$
      in the weak Poincar\'e inequality
      \beq\label{WP1} \mu(f^2)\le \aa(r)\EE(f,f)+ r \|f\|_\infty^2,\ \ \ f\in\D(\EE),\ \mu(f)=0.\end{equation}

\beg{thm}[\cite{RW01}, Theorems 2.1 and 2.3]\label{TT3}  If $(\ref{WP1})$ holds, then
$$ \|P_t-\mu\|_{\infty\to 2}^2  \le \xi(t),\ \ \ t\ge
0  $$
holds for $$\xi(t):=\inf\Big\{2r:\ r>0,\ -\ff 1 2 \aa(r)\log r\le t\Big\}$$
which goes to $0$    as $t\uparrow \infty$.

On the other hand, if $(\EE,\D(\EE))$ is symmetric and $$\xi(t):= \|P_t-\mu\|_{\infty\to 2}^2\to 0\ \text{as}\ t\to\infty,$$ then $(\ref{WP1})$  holds with
$$\aa(r) = 2r \inf_{s>0}\ff 1 s \xi^{-1}(s\exp[1-s/r]), \ \ r>0,$$
 where $\xi^{-1}(t):= \inf\{r> 0: \xi(r)\le t\}.$ \end{thm}

 When $(\EE,\D(\EE))$ is symmetric, we have  $\|P_t-\mu\|_{\infty\to 1}=\|P_t-\mu\|_{\infty\to 2}^2,$ so that
   the following is a consequence of Theorem \ref{TT3}.

 \beg{cor}[\cite{RW01}, Corollary 2.4] \label{CC3}  Let $(\EE,\D(\EE))$ be symmetric, then the following assertions hold.
\beg{enumerate} \item[$(1)$]   Let $\vv\in (0,1).$ Then $\eqref{WP1}$ holds with
 $$\aa(r)= c_1+c_2[\log(1+r^{-1})]^{(1-\vv)/\vv}$$ for some constants $c_1,c_2\in (0,\infty)$ if and only if
 $$\|P_t-\mu\|_{\infty\to 1}\le \exp[c_1'-c_2't^\vv]$$ holds for some constants $c_1',c_2'\in (0,\infty). $
 \item[$(2)$]   Let $p,q\in (1,\infty)$ with $p^{-1}+q^{-1}=1.$ Then $\eqref{WP1}$ holds with $$\aa(r)= c r^{1-p}$$ for some constant $c\in (0,\infty)$ if and only if
 $$\|P_t-\mu\|_{\infty\to 1}\le c' t^{1-q}$$ for some constant $c'\in (0,\infty)$.
 \item[$(3)$]   Let $p>0$. Then $\eqref{WP1}$ holds with
 $$\aa(r)=\exp[c(1+r^{-1/p})]$$
for some $c\in (0,\infty)$ if and only if
 $$\|P_t-\mu\|_{\infty\to 1}\le c'[\log(1+t)]^{-p}$$ holds for some   $c'\in (0,\infty).$\end{enumerate}
\end{cor}

\section{The case with boundary diffusion: $\dd>0$}

 When $\dd>0$, it is easy to estimate $\bb(r)$ by using $\bb_V(r)$ and $\bb_W^\pp(r).$ However, to estimate $\aa(r)$ using $\aa_V(r)$ and $\aa_W^\pp(r)$, we need the following assumption,
 where the function $h$ can be constructed by using the distance function to the boundary, see Example 3.1 and Example 4.1  for details.

 \beg{enumerate} \item[{\bf(A)}] There exists a function $h\in C^2(\bar M)$  such that
 $$Nh|_{\pp M}=1,\ \ \|\nn h\|_{L^2(\mu)}+ \|1_{\pp M} (W-V)^+\|_\infty+ \big\|( L_V h )^-\big\|_{L^2(\mu_V)}<\infty,$$
where $L_V:= \DD+\nn  V$ on $ M.$
\end{enumerate}
Under {\bf(A)}, we have
\beg{align*}&C_0  := \ff{Z_V\e^{\|1_{\pp M}(W-V)^+\|_\infty}} {Z_W^\pp }<\infty,\\
&C_1 :=\big\|(L_V h)^-\big\|_{L^2(\mu_V)} <\infty,\\
&C_2 := \|\nn h\|_{L^2(\mu_V)}<\infty.\end{align*}

\beg{thm}\label{T1} Let $\dd>0$.
\beg{enumerate}  \item[$(1)$] We have
\beq\label{E1}\bb(r)\le  \max\bigg\{\ff{\bb_V(r)}\theta,\  \ff{\bb_W^\pp(\dd r)}{1-\theta}\bigg\},\ \ r>0.\end{equation}
\item[$(2)$] If {\bf(A)} holds,    then  for any $r>0$,
\beq\beg{split} \label{E2}\aa(r)\le \inf\bigg\{& \max\Big\{  \big(1+ C_0^2C_1^2\big) \aa_V(s_1)+ \ff{1-\theta}\theta C_0^2C_1^2,\ \ff{1}\dd  \aa_W^\pp(s_2)\Big\}:\\
&\qquad  \ s_1,s_2>0,\   \theta\big(1 +C_0^2 C_1^2 \big)s_1+  (1-\theta)s_2 \le \ff r 4\bigg\}.\end{split}\end{equation}
In particular, taking
$$s_1=s_2:=s(r):=  \ff r{4 +4\theta C_0^2 C_1^2 },\ \ r>0,$$ we find a constant $c>1$ such that
\beq\label{E2'}\beg{split} \aa(r)&\le   \max\bigg\{  \big(1+ C_0^2C_1^2\big) \aa_V(s(r))+ \ff{1-\theta}\theta C_0^2C_2^2,\ \ff{1}\dd \aa_W^\pp(s(r))\bigg\}\\
&\le c+c \aa_V(r/c)+c\aa_W^\pp(r/c),\ \ \ r>0.\end{split} \end{equation}
\end{enumerate}
\end{thm}

\beg{proof} (1) For any $f\in C_b^1(\bar M)$, $\eqref{SP0}$ and \eqref{SP1'} imply
\beg{align*} &\mu(f^2)= \theta\mu_V(f^2)+(1-\theta)\mu_W^\pp(f^2)\\
&\le \theta r \mu_V(|\nn f|^2)+\theta \bb_V(r)\mu_V(|f|)^2+ (1-\theta)\dd r \mu_W^\pp(|\nn^\pp f|^2_\pp)
+(1-\theta) \bb_W^\pp(\dd r)\mu_W^\pp(|f|)^2\\
&\le r \big[\theta \mu_V(|\nn f|^2)+(1-\theta) \dd \mu_W^\pp(|\nn^\pp f|^2_\pp)\big] \\
 &\qquad +  \Big(\ff{\bb_V(r)}\theta\lor \ff{\bb_W^\pp(\dd r)}{1-\theta}\Big)\big[\theta\mu_V(|f|)+(1-\theta)\mu_W^\pp(|f|)\big]^2\\
&= r\EE_\dd(f,f)+ \Big(\ff{\bb_V(r)}\theta\lor \ff{\bb_W^\pp(\dd r)}{1-\theta}\Big)\mu(|f|)^2,\ \ r>0.\end{align*}
So, \eqref{E1} holds.

(2)  Let $f\in C_b^1(\bar M)$ such that
\beq\label{PO} \mu(f):= \theta\mu_V(f)+ (1-\theta)\mu_W^\pp(f)=0.\end{equation}   Then
\beg{align*}&\theta \mu_V(f)^2+(1-\theta)\mu_W^\pp(f)^2\\
&= \mu(f)^2 + \theta(1-\theta)|\mu_V(f)-\mu_W^\pp(f)|^2\\
&= \theta(1-\theta)|\mu_V(f)-\mu_W^\pp(f)|^2.\end{align*}
Combining this with \eqref{WP0} and \eqref{WP1'}, we derive
\beq\beg{split}\label{P1} &\mu(f^2)=\theta\mu_V(f^2)+(1-\theta)\mu_W^\pp(f^2)\\
&= \theta\mu_V(|f-\mu_V(f)|^2)+(1-\theta)\mu_W^\pp(|f-\mu_W^\pp(f)|^2)+\theta \mu_V(f)^2+(1-\theta)\mu_W^\pp(f)^2\\
&\le  \theta \aa_V(s_1) \mu_V(|\nn f|^2)+(1-\theta)\aa_W^\pp(s_2) \mu_W^\pp(|\nn^\pp f|_\pp^2)\\
&\quad + \big(\theta s_1+ (1-\theta)s_2\big)\big(\|f-\mu_V(f)\|_\infty^2\lor \|1_{\pp M}(f-\mu_W^\pp(f))\|_{\infty}^2\big) \\
&\quad +\theta(1-\theta)|\mu_V(f)-\mu_W^\pp(f)|^2.
\end{split}\end{equation} By $Nh|_{\pp M}=1$ and the integration by parts formula, we have
\beg{align*} &|\mu_V(f)-\mu_W^\pp(f)| = |\mu_W^\pp(f-\mu_V(f))| \\
&\le \ff{\e^{\|1_{\pp M}|(W-V)^+\|_\infty}}{Z_W^\pp}\int_{\pp M} |f-\mu_V(f)|\e^V\d\LL_\pp\\
&=  \ff{\e^{\|1_{\pp M}|(W-V)^+\|_\infty}}{Z_W^\pp}\int_{\pp M} |f-\mu_V(f)|  (Nh)\, \e^V\d\LL_\pp \\
&=-C_0 \int_M\Big[ |f-\mu_V(f)| L_V h +\<\nn  |f-\mu_V(f)|,\nn h\>\Big] \d\mu_V \\
&\le C_0 C_1 \|f-\mu_V(f)\|_{L^2(\mu_V)} +C_0C_2\|\nn f\|_{L^2(\mu_V)}.
\end{align*}
Noting that
$$\theta(1-\theta)(a+b)^2\le \theta a^2 +(1-\theta)b^2,\ \ a,b\ge 0,$$
this implies
\beq\label{GPP} \theta(1-\theta) |\mu_V(f)-\mu_W^\pp(f)|^2\le  \theta C_0^2C_1^2\mu_V(|f-\mu_V(f)|^2) + (1-\theta) C_0^2C_2^2 \mu_V(|\nn f|^2).\end{equation}
Combining this with \eqref{P1} and   $$\|f-\mu_V(f)\|_\infty\lor \|f-\mu_W^\pp(f)\|_\infty\le 2\|f\|_\infty,$$   for any $s_1,s_2>0$ and $f\in C_b^1(\bar M)$ with $\mu(f)=0$  we have
\beg{align*}&\mu(f)^2\le \theta \aa_V(s_1) \mu_V(|\nn f|^2)+(1-\theta)\aa_W^\pp(s_2) \mu_W^\pp(|\nn^\pp f|_\pp^2) \\
&\quad + 4\big(\theta s_1+ (1-\theta)s_2\big)\|f\|_\infty^2\\
&\quad +   \theta  C_0^2 C_1^2 \mu_V( |f-\mu_V(f)|^2) + (1-\theta) C_0^2 C_2^2  \mu_V(|\nn f|^2)\\
 &\le  \max\bigg\{  \big(1+ C_0^2C_1^2\big) \aa_V(s_1)+ \ff{1-\theta}\theta C_0^2C_1^2,\ \ff{1}\dd  \aa_W^\pp(s_2) \bigg\} \EE_\dd(f,f) \\
 &\quad +  4\Big[ \theta\big(1 +C_0^2 C_1^2 \big)s_1+  (1-\theta)s_2\Big]\|f\|_\infty^2.
 \end{align*}
Hence, \eqref{E2} holds, which implies \eqref{E2'} for the given choice of $s_1$ and $s_2$.
\end{proof}

To illustrate the above result, we present below an example to derive sharp functional inequalities for the sticky-reflected diffusion process, where the semigroup
is ultarbounded if and only if  $\tau>2,$ hypercontractive if and only if  $\tau=2$,   $L^2$-uniformly integrable if and only if  $\tau>1,$
$L^2$-exponential ergodic if and only if $\tau=1$, and sub-exponential ergodic when $\tau\in (0,1).$ For examples with weaker convergence rate, for instance the algebraic or logarithmic convergence,
one may take $V$ and $W$ with slower growth as in \cite[Example 1.4]{RW01}.

\paragraph{Example 3.1.}  Let $M:= (0,\infty)\times \R^{d}$ for some $d\ge 1$, and let $W(x)=V(x)= - |x|^\tau$ for some constant $ \tau>0.$ We have $\pp M=\{0\}\times\R^d\equiv \R^d.$

{\bf  (1)} It is known that \eqref{SP0} holds for some $\bb_V $ if and only if $\tau>1$, and in this case the exact order of $\bb_V(r)$ and $\bb_W^\pp(r)$ is
$\e^{cr^{-\ff\tau{2(\tau-1)}}}$  for small $r>0$, see  \cite[Corollary 2.5]{W00}. So, by Theorem \ref{T1}(1), there exists a constant $c_1>0$ such that the sticky-reflected diffusion process on
$\bar M:=[0,\infty)\times \R^d$ with $\pp M=\{0\}\times \R^d$ satisfies
$$\mu(f^2)\le r\EE_\dd(f,f)+ \exp\Big[c_1+c_1r^{-\ff{\tau}{2(\tau-1)}}\Big] \mu(|f|)^2,\ \ r>0,\ f\in C^1_b(\bar M),$$
that is,
\beq\label{BS} \bb(r)\le \exp\Big[c_1+c_1r^{-\ff{\tau}{2(\tau-1)}}\Big] \mu(|f|)^2,\ \ r>0.\end{equation}
This   is sharp for small $r>0$, since
it implies
$$\mu_V(f^2)\le r \mu_V(|\nn f|^2) + \exp\Big[c_1+c_1r^{-\ff{\tau}{2(\tau-1)}}\Big] \mu_V(|f|)^2,\ \ f\in C_0^1(M),$$
and this inequality  is sharp as shown in \cite[Corollary 2.5]{W00}.  Consequently:
\beg{enumerate}\item[(a)] If $\tau \in (1,2),$
then by Corollary \ref{CC1}, the Markov semigroup $P_t$ associated with $\EE_\dd$ is $L^2$-uniformly integrable with
$$\sup_{\mu(f^2)= 1} \int_{\bar M} (P_tf)^2 \exp\Big[C_t  \big\{\log (1+ (P_tf)^2)\big\}^{\ff{2(\tau-1)}\tau}\Big]\d\mu<\infty,\ \ t>0$$
for some $C: (0,\infty)\to (0,\infty).$
\item[(b)] If $\tau=2$, then  by \cite[Theorem 3.3.13(1)]{W05} the defective log-Sobolev inequality holds, which together with   \cite[Corollary 1.2]{W14} implies the strict log-Sobolev inequality,
since it is well known that the defective log-Sobolev inequality together with the Poincar\'e  inequality implies the strict log-Sobolev inequality. So, according to \cite{Gross1} or \cite{Gross2},  in this case
 $P_t$ is hypercontractive, i.e.  $\|P_t\|_{L^2(\mu)\to L^4(\mu)]}\le 1$ holds  for large enough $t>0$.
  \item[(c)] If $\tau>2$, then by Corollary \ref{CC2}, $P_t$ is ultrabounded with
$$\|P_t\|_{L^1(\mu)\to L^\infty(\mu)}\le \e^{c + ct^{-\ff\tau{\tau-2}}},\ \ \ t>0$$
  for some constant $c>0.$ \end{enumerate}

{\bf (2)}  Next, $\mu_V$ and $\mu_W^\pp$ satisfy the Poincar\'e inequality if and only if $\tau=1$, see \cite{W99},  so that in this case $\EE_\dd$ satisfies the Poincar\'e inequality as well, due to
Theorem \ref{T1}(2) with $h(r,x):= r$ for $(r,x)\in \bar M$ and bounded $\aa_V$ and $\aa_W^\pp.$  Consequently,  when $\tau=1$,
$$\|P_t-\mu\|_{L^2(\mu)}\le \e^{-\ll t},\ \ \ t\ge 0$$ holds for some constant $\ll>0.$

{\bf (3)}  Moreover, when $\tau\in (0,1)$, by \cite[Example 1.4(c)]{RW01} which applies also to $(0,\infty)\times \R^d$ in place of $\R^{d+1}$, $\aa_V(r)$ and $\aa_W^\pp(r)$ behaves as $[\log r^{-1}]^{\ff{4(1-\tau)}\tau}$ for small $r>0$,
so that by Theorem \ref{T1}(2) for $h(r,x):= r$ for $ (r,x)\in \bar M$,
there exists a constant $c_2>0$ such that
$$\mu(f^2)\le [\log (c_2+r^{-1})]^{\ff{4(1-\tau)}\tau} \EE_\dd(f,f)+ r\|f\|_\infty^2,\ \ r>0,\ f\in C^1_b(\bar M),\ \mu(f)=0.$$
Consequently, by Corollary \ref{CC3}, the Markov semigroup $P_t$ associated with $\EE_\dd$ is sub-exponential ergodic
$$\|P_t-\mu\|_{L^\infty(\mu)\to L^2(\mu)} \le \e^{-c t^{\ff\tau{4-3\tau}}},\ \ t>0$$
for some constant $c>0$.

\section{The case without boundary diffusion: $\dd=0$}

When $\dd=0$, the Dirichlet form for the sticky-reflected diffusion reduces to
$$\EE_\dd(f,f)=\mathscr E_{0} (f,f):=\mu(|\nn f|^2) =\theta\mu_V(|\nn f|^2),\ \ f\in C^1_b(\bar M).$$
So, to establish the weak and super Poincar\'e inequalities, we need to bound the $L^2(\mu_W^\pp)$-norm using the Neumann Dirichelt form $\mu_V(|\nn f|^2).$
To this end, we need the following assumption, which is slightly stronger than {\bf(A)}.

 \beg{enumerate} \item[{\bf (B)}] There exists a function $h\in C^2(\bar M)$  such that
 $$Nh|_{\pp M}=1,\ \ \|\nn h\|_{\infty)}+ \|1_{\pp M} (W-V)^+\|_\infty+ \big\|( L_V h )^-\big\|_{\infty}<\infty,$$
where $L_V:= \DD+\nn  V$ on $ M$.
\end{enumerate}
Under {\bf (B)}, besides
$C_0  := \ff{Z_V\e^{\|1_{\pp M}(W-V)^+\|_\infty}} {Z_W^\pp }<\infty$, we have
$$\bar C_1 :=\big\|(L_V h)^-\big\|_{\infty} <\infty,\ \ \ \
 \bar C_2:= \|\nn h\|_{\infty}<\infty.$$

\beg{thm}\label{T2} Assume {\bf(B)} and let $\dd=0. $  Then
\beq\label{B1} \bb(r)\le   \Big(\ff 2{\theta r} C_0^2\bar C_2^2 + C_0\bar C_1\Big)\bb_V\bigg(\ff{\theta^2r^2}{4C_0^2\bar C_2^2+ 2\theta C_0\bar C_1 r}\bigg),\ \ r>0,\end{equation}
 \beq\label{B2} \aa(r)\le   \ff{A}\theta\aa_V\Big(\ff{r}{4A}\Big)+ \ff B\theta,\end{equation}
 where
 \beg{align*}&A:= \theta+(1-\theta) \big(C_0^2\bar C_2^2+ C_0\bar C_1\big)\big((1-\theta)C_0^2C_1^2+1\big),\\
 &B:= 1-\theta+\ff{(1-\theta)^2}\theta C_0^2C_2^2.\end{align*}

\end{thm}

\beg{proof} (a) By $Nh|_{\pp M}=1$, and the integration by parts formula,
we obtain
\beq\label{PP}\beg{split}  \mu_W^\pp (f^2)&= \mu_W^\pp (f^2 Nh)\le \ff{\e^{\|1_{\pp M}(W-V)^+ \|_\infty} }{Z_W^\pp} \int_{\pp M} f^2( Nh )\e^{V}\d\LL_\pp\\
&= - C_0 \mu_V \big(f^2  L_V h+\<\nn f^2,\nn h\>\big)\\
&\le C_0\bar C_1  \mu_V(f^2) + 2C_0\bar C_2 \ss{\mu_V(f^2) \mu_V(|\nn f|^2)}\\
&\le s \mu_V(|\nn f|^2) + \big(s^{-1} C_0^2\bar C_2^2 + C_0\bar C_1  \big) \mu_V(f^2),\ \ s>0.\end{split}\end{equation}
Combining this with \eqref{SP0}, we obtain
\beg{align*}  \mu_W^\pp (f^2)  &\le 2s  \mu_V(|\nn f|^2) +\big(s^{-1} C_0^2 \bar C_2^2+ C_0\bar C_1\big) \bb_V\Big(\ff s{s^{-1} C_0^2\bar C_2^2 +C_0\bar C_1}\Big)\mu_V(|f|)^2\\
&= \ff{2s}\theta  \EE_0(f,f)  +\big(s^{-1} C_0^2 \bar C_2^2+ C_0\bar C_1\big) \bb_V\Big(\ff {s^2}{C_0^2\bar C_2^2 +sC_0\bar C_1}\Big)\mu_V(|f|)^2.\end{align*}
By taking $s=\ff \theta 2 r,$  we  derive    the   estimate \eqref{B1}  on $\bb(r)$.

(b) Let $f\in C_b^1(\bar M)$ with $\mu(f)=\theta\mu_V(f)+(1-\theta)\mu_W^\pp (f)=0.$ Then by the triangle inequality,
\beg{align*} |\mu_V(f)|&\le |\theta\mu_V(f)+ (1-\theta)\mu_W^\pp(f)| +|(1-\theta)\mu_V(f)- (1-\theta) \mu_W^\pp(f)|\\
&= (1-\theta)|\mu_V(f)-\mu_W^\pp(f)|.\end{align*}
Combining this with \eqref{GPP} we derive
$$\mu_V(f)^2 \le (1-\theta)C_0^2C_1^2\mu_V(|f-\mu_V(f)|^2) + \ff{(1-\theta)^2}\theta C_0^2C_2^2\mu_V(|\nn f|^2).$$
Combining this with \eqref{WP0} and
$$\mu_V(f^2)=\mu_V(|f-\mu_V(f)|^2)+\mu_V(f)^2,$$
we obtain
\beq\label{OI} \beg{split}&\mu_V(f^2)  \le \big(1+ (1-\theta)C_0^2C_1^2\big)\mu_V(|f-\mu_V(f)|^2) + \ff{(1-\theta)^2}\theta C_0^2C_2^2\mu_V(|\nn f|^2)\\
&\le \Big[\ff{(1-\theta)^2}\theta C_0^2C_2^2+  \big(1+ (1-\theta)C_0^2C_1^2\big)\aa_V(s)\Big] \mu_V(|\nn f|^2)\\
&\qquad + \big(1+ (1-\theta)C_0^2C_1^2\big) s \|f-\mu_V(f)\|_\infty^2,\ \ s>0.\end{split}
\end{equation}
On the other hand, taking $s=1$ in \eqref{PP}, we obtain
\beg{align*} \mu(f^2)&=\theta\mu_V(f^2)+(1-\theta)\mu_W^\pp(f^2)\\
&\le (1-\theta) \mu_V(|\nn f|^2)+ \big[\theta+(1-\theta) (C_0\bar C_2^2+C_0\bar C_1)\big]\mu_V(f^2).\end{align*}
Combining this with   \eqref{OI}, $\|f-\mu_V(f)\|_\infty\le 2\|f\|_\infty$,  and the definitions of $A$ and $B$, we arrive at
\beg{align*}\mu(f^2)&\le \big(B+A\aa_V(s)\big)\mu_V(|\nn f|^2) + As\|f-\mu_V(f)\|_\infty^2\\
&\le \Big(\ff B\theta +\ff A\theta \aa_V(s)\Big)\mathscr E_0(f,f)+ 4 As \|f\|_\infty^2,\ \ s>0.\end{align*}
Taking $s=\ff r{4A}$ we obtain \eqref{B2}.
\end{proof}

By Theorem \ref{T2}, for the model in Example 3.1 with $\dd=0$,  the same assertion for $\tau\in (0,1]$ holds as in Example 3.1, and when $\tau>1$ we have
$$\bb(r)\le \exp\Big[c_1+c_1r^{-\ff\tau{\tau-1}}\Big],\ \ \ r>0$$
for some constant $c_1>0$,
which is weaker than the corresponding ones for $\dd>0$.

Below we consider the sticky-reflected diffusion on   a compact manifold.

\paragraph{Example 4.1.} Let $\bar M$ be a d-dimensional compact Riemannian manifold with smooth boundary $\pp M$, and let $\dd=0$. Then
$\eqref{SP}$ holds for
\beq\label{X0} \bb(r)=  c (1\land r)^{-d},\ \ r>0\end{equation}
for some constant $c>0$.  When $d=1$ \eqref{X0} can be improved as
\beq\label{X0'}  \bb(r)=  c (1\land r)^{-\ff 1 2}= \bb(r)=  c (1\land r)^{-\ff d2}.\end{equation}

\beg{proof} (a) Since $W$ and $V$ are bounded as $\bar M$ is compact, it suffices to consider $W=V=0.$ So,
$$\mu_V(\d x)=\ff{\LL(\d x)}{\LL(M)}=:\mu_0(\d x),\ \ \mu_W^\pp(\d x)=\ff{\LL_\pp(\d x)}{\LL_\pp(\pp M)}=:\mu_0^\pp(\d x).$$
Let $\rr_\pp$ be the Riemannian distance to the boundary $\pp M$. Then there exists a constant $s_0>0$ such that   $\rr_\pp\in C_b^2(\pp_{s_0}M),$
where
$$\pp_{s_0}M:= \{x\in \bar M:\ \rr_\pp(x)\le s_0\}.$$
Consider the polar coordinates on $\pp_{s_0}M$
$$[0,s_0]\times\pp M\ni (r,z)\mapsto \exp_z[r N(z)]\in  \pp_{s_0}M.$$
Then there exists a constant $c_0>1$ such that the volume measure $\LL$ on $\pp_{s_0}M$ satisfies
\beq\label{CD} c_0^{-1}  \d r \LL_\pp(\d z)  \le \LL(\d r,\d z)\le c_0 \d r \LL_\pp(\d z).\end{equation}
Choose   $\xi\in C^\infty([0,\infty); [0,1])$ such that $\xi'\ge 0, \xi(0)=1$ and $\xi(s)=0$ for $s\ge s_0$, and let
$$h (x):= \beg{cases} r \xi (r),\ &\text{if}\ x=\exp_z(rN(z))\in \pp_{s_0}M,\\
0,  &\text{if}\ x\in M\setminus \pp_{s_0}M.\end{cases} $$
Then $h \in C_b^2(\bar M)$ with $h|_{\pp M}=0$ and $Nh|_{\pp M}=1.$  So, the assumption {\bf (B)} holds.

Now, since $M$ is a $d$-dimensional compact manifold, the classical Nash inequality with dimension $d$ holds for the Dirichlet form $\mu_V(|\nn f|^2)$, so by \cite[Corollary 3.3]{W00},
there exists a constant $c_1>0$ such that
$$\bb_V (r)\le c_1 (1\land r^{-\ff d 2}),\ \ \ r>0.$$
Then the desired assertion follows from Theorem \ref{T2}.

(b) When $d=1$, we may simply consider $M=(0,1)$ and $\gg=\ff 1 2$ so that $\theta=\ff 12$, hence
$$\mu(f)= \ff 1 2\int_0^1f(s)\d s+ \ff 1 4f(0)+\ff 1 4 f(1)$$ and
$$\EE_0(f,f)= \ff 1 2 \int_0^1 f'(s)^2\d s.$$
 By the classical Nash inequality, there exists a constant $c_1>0$ such that
  $$\int_0^1 f(s)^2\d s\le r \int_0^1 f'(s)^2\d s+ c_1 (1\land r)^{-\ff 1 2} \bigg(\int_0^1|f(s)|\d s\bigg)^2,\ \ s>0,\ f\in C_b^1([0.1]).$$
Then there exists a constant $c>0$ such that
\beg{align*} &\mu(f^2)= \ff 1 2\int_0^1 f(s)^2\d s+ \ff 1 4 f(0)^2+\ff 1 4 f(1)^2\\
&\le \ff  r 2\int_0^1f'(s)^2\d s+ \ff 1 2 c_1 (1\land r)^{-\ff 1 2} + \ff 1 4 f(0)^2+\ff 1 4 f(1)^2\\
&\le r \EE_0(f,f) + c (1\land r)^{-\ff 1 2} \mu(|f|)^2,\ \ \ r>0,\ f\in C_b^1([0,1]).\end{align*}
Then \eqref{SP} holds for $\dd=0$ and $\bb(r)$ in \eqref{X0'}.
\end{proof}

 \paragraph{Problem 4.1.} We hope that in Example 4.1, even for $d\ge 2$ we still have \eqref{SP}  for $\bb$ in \eqref{X0'}. To this end, the general estimate \eqref{B1}  should be improved by more  refined calculus.
\ \newline
 {\bf Acknowledgement.} The work is supported in part by  National Key R\&D Program of China (No. 2022YFA1006000) and NSFC (12531007). The author would like to thank the referee for helpful comments.

\end{document}